\documentclass[11pt]{article}
\usepackage{amssymb,amsthm,amsmath,amsthm}
\usepackage[numbers,sort&compress]{natbib}
\newtheorem{theorem.}{\textbf{Theorem}}[section]
\newtheorem{lemma.}{\textbf{Lemma}}[section]
\newtheorem{remark.}{\textbf{Remark}}
\date{}
\begin{document}

\title{\Large{\bf New homoclinic orbits for Hamiltonian systems with asymptotically quadratic growth at infinity}\thanks{D.-L. Wu is supported by NSF of China (No.11801472), China Scholarship Council (No.201708515186)
and the Youth Science and Technology Innovation Team of Southwest
Petroleum University for Nonlinear Systems (No.2017CXTD02). X. Yu is
supported by NSF of China (No.11701464) and the Fundamental Research
Funds for the Central Universities (No.JBK1805001).}}

\author{Dong-Lun Wu$^{a,b}$, Xiang Yu$^{c}$\footnote{Corresponding author.}\\
{\small $^{a}$College of Science, Southwest Petroleum University,}\\
{\small Chengdu, Sichuan 610500, P.R. China}\\
{\small $^{b}$Institute of Nonlinear Dynamics, Southwest Petroleum University,}\\
{\small Chengdu, Sichuan 610500, P.R. China}\\
{\small $^{c}$School of Economic and Mathematics, Southwestern
University of Finance and Economics,}\\
{\small Chengdu, Sichuan 611130, P.R. China}\\
}
\date{}
\maketitle

\let\thefootnote\relax\footnotetext{Email address: wudl2008@163.com; yuxiang@swufe.edu.cn}

{\bf Abstract} In this paper, we study the existence and multiplicity
of homoclinic solutions for following Hamiltonian systems with
asymptotically quadratic nonlinearities at infinity
\begin{eqnarray*}
\ddot{u}(t)-L(t)u+\nabla W(t,u)=0.
\end{eqnarray*}
We introduce a new coercive condition and obtain a new embedding
theorem. With this theorem, we show that above systems possess at
least one nontrivial homoclinic orbits by Generalized Mountain Pass
Theorem. By Variant Fountain Theorem, infinitely many homoclinic
orbits are obtained for above problem with symmetric condition. Our
asymptotically quadratic conditions are different from
 previous ones in the references.


\section{Introduction}

In this paper, we consider the following systems
\begin{eqnarray}
\ddot{u}(t)-L(t)u+\nabla W(t,u)=0\label{1},
\end{eqnarray}
where $W\in C^{1}(\mathbb{R}\times\mathbb{R}^{N},\mathbb{R})$ and
$L\in C(\mathbb{R},\mathbb{R}^{N^{2}})$ is a symmetric matrix. A
solution $u(t)$ of (\ref{1}) is called nontrivial homoclinic if
$u\not\equiv0$, $u(t)\rightarrow0$ as $t\rightarrow\pm\infty$.
Moreover, given a $N\times N$ matrix $\mathcal{M}$, we say
 $\mathcal{M}\geq0$ if $\inf_{|x|=1}(\mathcal{M}x,x)\geq0$ and
$\mathcal{M}\ngeq0$ if $\mathcal{M}\geq0$ does not hold.

In last decades, along with the development of variational methods,
many mathematicians showed the existence and multiplicity of solutions
for differential
equations(see\cite{20,52,1,2,3,53,43,50,42,17,21,6,41,7,19,40,8,10,11,31,13,60,61,51,14,18,15,16,30,70,71,72,73,74,75,76,77}).
One of the difficulties to obtain homoclinic orbits for (\ref{1}) is
 the lack of compactness of embeddings. To solve this problem,
Rabinowitz and his co-authors introduced the periodic and coercive
conditions. With the periodic assumption, Rabinowitz \cite{41}
obtained homoclinic solutions for (\ref{1}) by taking limit of a
sequence of subharmonic orbits. After then there are many papers
concerning on the existence and multiplicity of homoclinic solutions
for problem (\ref{1}) under periodic assumption. Without periodic
assumptions, Rabinowitz and Tanaka \cite{7} introduced following
coercive condition on $L(t)$ to retrieve the compactness.

\vspace{0.2cm} $(L')$ Let $L$ be a positively definite matrix and
\begin{eqnarray*}
l(t)=\inf_{|x|=1}(L(t)x,x)\rightarrow+\infty\
\ \mbox{as}\ \ |t|\rightarrow\infty.
\end{eqnarray*}
Omana and Willem \cite{42} used condition $(L')$ to obtain a compact
embedding theorem. In 2010, Wang, Zhang and Xu \cite{14} introduced
the following new coercive condition.

\vspace{0.2cm}$(L'')$ The measure of the set
$\Lambda^{b}:=\{t\in\mathbb{R}: L(t)<b I_{N}\}$ is finite for some
$b>0$, where $I_{N}$ is the identity matrix of order $N$.

Under condition $(L'')$, the spectrum of the operator
$-(d^{2}/dt^{2})+L(t)$ can go to $-\infty$ and the corresponding
functional of problem (\ref{1}) becomes strong indefinite, i.e.
unbounded from below and from above on infinite dimensional spaces.
Another coercive condition can be found in \cite{43} as follow.

\vspace{0.2cm}$(L''')$ Form some $r_{0}>0$,
\begin{eqnarray*}
\lim_{|s|\rightarrow+\infty}\mbox{meas}\left(\left\{t\in(s-r_{0},s+r_{0}):L(t)\ngeq y I_{N}\right\}\right)=0,\ \ \ \forall y>0.
\end{eqnarray*}

In 1995, in order to get the
 homoclinic solutions for subquadratic
Hamiltonian systems, Ding \cite{2} introduced the following condition.

\vspace{0.2cm}$(L'''')$ From some $\varsigma<1$,
$l(t)|t|^{\varsigma-2}\rightarrow+\infty$ as $|t|\rightarrow\infty$.

\vspace{0.2cm} As far as we know, in \cite{2}, Ding has considered the
case $L(t)$ is not positively definite for the first time. Combining
with the coercive condition, $L(t)$ possesses finite negative
eigenvalues. Using similar conditions, in \cite{1}, Chen obtained the
ground state homoclinic solution for (\ref{1}) when 0 lies in the gap
of spectrum of $-(d^{2}/dt^{2})+L(t)$. Recently, Schechter \cite{19}
studied problem (\ref{1}) with the conditions weak enough to make the
linear operator have the essential spectrum possibly and the author
considered different cases according to the number of the negative
eigenvalues.

As we know, the growth of nonlinear term $W$ is important in showing
the geometric structure of the corresponding functional and the
boundedness of the asymptotic critical points sequence. The growth of
$W$ is mainly classified into
 superquadratic, subquadratic and asymptotically quadratic cases. In this paper, we consider problem (\ref{1}) with a set of new
asymptotically quadratic growth conditions at infinity. There are
already many results concerning on homoclinic solutions for the
Hamiltonian systems with asymptotically quadratic growth. By
investigating the references, we can find that, when dealing with the
asymptotically quadratic Hamiltonian systems, some mathematicians
considered the following condition.

\vspace{0.2cm}$(G1)$ $\widetilde{W}(t,x)<0$ if $x\neq0$ and
$\widetilde{W}(t,x)\rightarrow-\infty$ as $|x|\rightarrow\infty$
uniformly in $t$, and
\begin{eqnarray*}
\inf_{t\in\mathbb{R}, c_{1}\leq|x|\leq c_{2}}\frac{\widetilde{W}(t,x)}{|x|^{2}}<0
\end{eqnarray*}
for all $0<c_{1}<c_{2}<+\infty$, where
$\widetilde{W}(t,x)=2W(t,x)-(\nabla W(t,x),x)$.

Let $\mathcal{H}$ be the self-adjoint extension of
$-(d^{2}/dt^{2})+L(t)$ with $\mathfrak{D}(\mathcal{H})\subset
L^{2}(\mathbb{R},\mathbb{R}^{N})$. Some other mathematicians \cite{50}
assumed that

$(G2)$ Let $W$ satisfy
\begin{eqnarray*}
\nabla W(t,x)=a(t)x+\nabla G(t,x),
\end{eqnarray*}
where  $a(t)$ satisfies
\begin{eqnarray*}
\inf_{t\in\mathbb{R}}\left[\inf_{|x|=1}(a(t)x,x)\right]>\inf(\sigma(\mathcal{H})\cap(0,+\infty)),
\end{eqnarray*}
where $\sigma(\mathcal{H})$ the spectrum of $\mathcal{H}$ and $G(t,x)$
is subquadratic at infinity in $x$.

By above conditions and some other auxiliary conditions,
mathematicians obtained the homoclinic solutions for (\ref{1}). More
details can be found in \cite{52,2,3,50,19,40,8,13,51,15,16}. However,
in these papers, the authors required that

\vspace{0.2cm}$(G3)$ $\widetilde{W}(t,x)\leq0$ for all
$t\in\mathbb{R}$ with $|x|$ large enough.

\vspace{0.2cm}In this paper, we consider the nonlinearities with
asymptotic quadratic growth at infinity with $(G3)$ being not
satisfied. To regain the compactness, we introduce the following new
coercive condition.

$(L1)$ there exist $\alpha>1$ and $M>0$ such that
\begin{eqnarray*}
\lim_{R\rightarrow+\infty}\mbox{meas}\left(\left\{|t|\geq R:\frac{L(t)}{|t|\ln^{\alpha}|t|}\leq M I_{N}\right\}\right)=0.
\end{eqnarray*}

Besides, we assume the following condition.

$(L2)$ there exists $L_{0}>0$ such that $(L(t)x,x)\geq-L_{0}|x|^{2}$
for all $(t,x)\in\mathbb{R}\times\mathbb{R}^{N}$.

By $(L1)$ and $(L2)$, similar to \cite{2}, we can set
$E=\mathfrak{D}(|\mathcal{H}|^{1/2})$ be the domain of
$|\mathcal{H}|^{1/2}$ and the inner product and norm on $E$ are
$(u,w)_{o}=(|\mathcal{H}|^{1/2}u,|\mathcal{H}|^{1/2}w)_{2}+(u,w)_{2}$,
 which shows that $E$ is a Hilbert space. It
can be proved that $\sigma(\mathcal{H})$ consists of a sequence of
eigenvalues $\lambda_{1}\leq\lambda_{2} \leq \cdots\rightarrow\infty$,
and a sequence of eigenfunctions
$(e_{n})$$(\mathcal{H}e_{n}=\lambda_{n}e_{n})$ forms an orthogonal
basis in $L^{2}$, where  $\sigma(\mathcal{H})$ is the spectrum of
$\mathcal{H}$. Let $n^{-}$, $n^{0}$ and $\bar{n}$ stand for the
numbers of the negative, null and nonpositive eigenvalues
respectively. Set $E^{-} $, $E^{0}$ and $E^{+}$ be the negative, null
and positive space spanned with negative, null and positive
eigenvectors respectively.  Then $E=E^{-}\oplus E^{0}\oplus E^{+}$. We
can define
\begin{eqnarray*}
(u,w)=(|\mathcal{H}|^{1/2}u,|\mathcal{H}|^{1/2}w)_{2}+(u^{0},w^{0})_{2}\ \ \ \mbox{and}\ \ \ \|u\|=(u,u)^{1/2},
\end{eqnarray*}
where $u=u^{-}+u^{0}+u^{+}$ and $w=w^{-}+w^{0}+w^{+}\in E$. Clearly,
$\|\cdot\|_{o}$ is equivalent to $\|\cdot\|$ according to the proof of
Lemma 2.1 in \cite{2}. Our main results are as follow.

\begin{theorem.}\label{th1.1} Let $(L1)$, $(L2)$ hold and $W(t,x)$
satisfy

$(W1)$ $W(t,0)=0$, $\forall t\in\mathbb{R}$, there exists $b_{1}>0$
such that
\begin{eqnarray*}
|\nabla W(t,x)|\leq b_{1}|x|\ \ \ \forall (t,x)\in\mathbb{R}\times\mathbb{R}^{N};
\end{eqnarray*}

$(W2)$ there exist $b_{2}$, $r_{1}>0$ and $\mu>1$ such that
\begin{eqnarray*}
\widetilde{W}(t,x)\geq b_{2}|x|^{\mu}\ \ \mbox{if}\ \ |x|\geq r_{1}\ \ \ \forall t\in\mathbb{R};
\end{eqnarray*}

$(W3)$ there exists $b_{3}>0$ such that
\begin{eqnarray*}
\widetilde{W}(t,x)\geq -b_{3}|x|\ \ \ \forall t\in\mathbb{R}\ \ \mbox{and}\ \ |x|< r_{1};
\end{eqnarray*}

$(W4)$ there exist $r_{2}$, $\sigma_{0}>0$ such that
\begin{eqnarray*}
|W(t,x)|\leq\frac{\lambda_{\bar{n}+1}-\sigma_{0}}{2}|x|^{2}\ \ \   \mbox{for any}\ \
|x|\leq r_{2}\ \ \mbox{and}\ \ t\in\mathbb{R};
\end{eqnarray*}

$(W5)$ there exist $\varepsilon_{0}>0$, $r_{\infty}>0$ such that
\begin{eqnarray*}
W(t,x)\geq \left(\frac{\lambda_{\bar{n}+1}}{2}+\varepsilon_{0}\right)|x|^{2}\ \ \ \mbox{for any}\ \ |x|\geq r_{\infty}\ \ \mbox{and}\ \ t\in\mathbb{R};
\end{eqnarray*}

$(W6)$ $W(t,x)\geq0$ for all $(t,x)\in\mathbb{R}\times\mathbb{R}^{N}.$

Then there exists a homoclinic solution for (\ref{1}).
\end{theorem.}

\begin{theorem.}\label{th1.2}
Let $(L1)$-$(L2)$, $(W1)$-$(W6)$ hold and $W(t,x)$ is even in $x$,
then there are infinitely many homoclinic solutions for (\ref{1}).
\end{theorem.}

\begin{remark.}
Compared to condition $(L'''')$, $l(t)$ may not have limit under
conditions $(L1)$ and $(L2)$. Moreover, our conditions are also
different from $(L')$, $(L'')$ and $(L''')$.
\end{remark.}
\begin{remark.}
In our theorems, we introduce a set of new asymptotically quadratic
growth conditions where the potentials $W(t,x)$ does not satisfy
 $(G3)$.
\end{remark.}
\begin{remark.}
As we know, there are many papers in which the infinitely many
homoclinic solutions are obtained for superquadratic and subquadratic
Hamiltonian systems. However, there are only few similar results for
(\ref{1}) under asymptotically quadratic growth conditions. In Theorem
1.2, we conclude that problem (\ref{1}) has infinitely many homoclinic
solutions under asymptotically quadratic growth condition without
periodic conditions.
\end{remark.}

\section{Preliminaries}

Similar to \cite{2}, let $I: E \rightarrow \mathbb{R}$ be as follow.
\begin{eqnarray}
I(u)&=&\int_{\mathbb{R}}\left(\frac{1}{2}|\dot{u}(t)|^{2}+\frac{1}{2}(L(t)u(t),u(t))-W(t,u(t))\right)dt{\nonumber}\\
&=&
\mbox{}\frac{1}{2}\|u^{+}\|^{2}-\frac{1}{2}\|u^{-}\|^{2}-\int_{\mathbb{R}}W(t,u(t))dt.\label{7}
\end{eqnarray}

\begin{lemma.}\label{le2.1}Assume $(L1)$ and $(L2)$ hold, then $E\hookrightarrow L^{p}(\mathbb{R},\mathbb{R}^{N})$ is
compact for any $p\in[1,+\infty]$ and there is $C_{p}>0$ such that
\begin{eqnarray}
\ \|u\|_{L^{p}}\leq C_{p}\|u\|,\ \ \  \ \forall u\in
E.\label{4}
\end{eqnarray}
\end{lemma.}

{\bf Proof.} We adopt some ideas from \cite{20} and \cite{2}. For
$R>0$, set
\begin{eqnarray*}
\Gamma_{R}=\left\{t\in \mathbb{R}:|t|\leq R\right\}
\end{eqnarray*}
and
\begin{eqnarray*}
D_{M}(R)=\left\{t\in\Gamma_{R}^{c}:\frac{L(t)}{|t|\ln^{\alpha}|t|}\leq M
I_{N}\right\}.
\end{eqnarray*}
Since $\|\cdot\|_{o}$ is equivalent to $\|\cdot\|$, we use the norm
$\|\cdot\|_{o}$ to show the proof. First, we suppose that $l(t)\geq1$
for all $t\in\mathbb{R}$. Then $|\mathcal{H}|=\mathcal{H}$ and
\begin{eqnarray*}
\|u\|_{o}^{2}=\int_{\mathbb{R}}\left(|\dot{u}(t)|^{2}+(L(t)u(t),u(t))+|u(t)|^{2}\right)dt.
\end{eqnarray*}
Let $K\subset E$ satisfy $\|u\|\leq M_{1}$ for all $u\in K$ for some
$M_{1}>0$. Next, we show $K$ is precompact in
$L^{p}(\mathbb{R},\mathbb{R}^{N})$ for any $p\in[1,+\infty]$.

First, we put $p=1$. By Soblolev embedding compact theorem, for any
$\varepsilon>0$, there are $u_{1}$, $\ldots$, $u_{m}\in K$ such that
for $\forall u\in K$, we can find $u_{i}$($i\in[1,m]$) such that
\begin{eqnarray*}
\int_{\Gamma_{R}}|v_{i}(t)|dt\leq\frac{1}{2}\varepsilon,
\end{eqnarray*}
where $v_{i}(t)=u(t)-u_{i}(t)$. Obviously, one has
\begin{eqnarray}
\int_{\Gamma_{R}^{c}}|v_{i}(t)|dt\leq\int_{D_{M}(R)}|v_{i}(t)|dt+\int_{\Gamma_{R}^{c}\setminus D_{M}(R)}|v_{i}(t)|dt.\label{16}
\end{eqnarray}
Letting $1<\eta<\alpha$ and $\Omega_{n}=\{t\in\mathbb{R}:
(|t|\ln^{\eta}|t|)|v_{i}|\leq1\}$, we obtain that
\begin{eqnarray}
&&\int_{\Gamma_{R}^{c}\setminus D_{M}(R)}|v_{i}(t)|dt\nonumber\\
&\leq&\int_{(\Gamma_{R}^{c}\setminus D_{M}(R))\bigcap\Omega_{n}}\frac{1}{|t|\ln^{\eta}|t|}dt+\int_{(\Gamma_{R}^{c}\setminus D_{M}(R))\bigcap \Omega_{n}^{c}}(|t|\ln^{\eta}|t|)|v_{i}(t)|(|t|\ln^{\eta}|t|)^{-1}dt\nonumber\\
&\leq&\int_{(\Gamma_{R}^{c}\setminus D_{M}(R))\bigcap\Omega_{n}}\frac{1}{|t|\ln^{\eta}|t|}dt+\int_{(\Gamma_{R}^{c}\setminus D_{M}(R))\bigcap \Omega_{n}^{c}}(\ln^{\eta-\alpha}|t|)(|t|\ln^{\alpha}|t|)|v_{i}(t)|^{2}dt\nonumber\\
&\leq&\int_{\Gamma_{R}^{c}\setminus D_{M}(R)}\frac{1}{|t|\ln^{\eta}|t|}dt+\frac{1}{M}\int_{\Gamma_{R}^{c}\setminus D_{M}(R)}(\ln^{\eta-\alpha}|t|)(L(t)v_{i}(t),v_{i}(t))dt.\label{18}
\end{eqnarray}
From (\ref{16}) and (\ref{18}), there is $M_{2}>0$ such that
\begin{eqnarray*}
\int_{\Gamma_{R}^{c}}|v_{i}(t)|dt
&\leq&\int_{D_{M}(R)}|v_{i}(t)|dt+\int_{\Gamma_{R}^{c}\setminus D_{M}(R)}|v_{i}(t)|dt\nonumber\\
&\leq&\sqrt{\mbox{meas}(D_{M}(R))}\left(\int_{D_{M}(R)}|v_{i}(t)|^{2}dt\right)^{\frac{1}{2}}+\int_{\Gamma_{R}^{c}\setminus D_{M}(R)}\frac{1}{|t|\ln^{\eta}|t|}dt\nonumber\\
&&+\frac{1}{M}\int_{\Gamma_{R}^{c}\setminus D_{M}(R)}(\ln^{\eta-\alpha}|t|)(L(t)v_{i}(t),v_{i}(t))dt\nonumber\\
&\leq&\sqrt{\mbox{meas}(D_{M}(R))}\|v_{i}\|_{o}+\int_{\Gamma_{R}^{c}}\frac{1}{|t|\ln^{\eta}|t|}dt+\frac{1}{M\ln^{\alpha-\eta}R}\|v_{i}\|_{o}^{2}\nonumber\\
&\leq&2M_{1}\sqrt{\mbox{meas}(D_{M}(R))}+\frac{M_{2}}{\ln^{\eta-1}R}+\frac{M_{1}}{M\ln^{\alpha-\eta}R}.
\end{eqnarray*}
From $(L1)$, we can choose $R$ large enough such that $\int_{\{t\in
\mathbb{R}: |t|\geq R\}}|v_{i}(t)|dt<\frac{1}{2}\varepsilon$. Then we
can see $K$ has a finite $\varepsilon$-net. Therefore, $K$ is
precompact in $L^{1}(\mathbb{R},\mathbb{R}^{N})$.

Second, we put $p=+\infty$. For $\forall n\in \mathbb{N}$ and
$t\in\mathbb{R}$, we have
\begin{eqnarray*}
u(t)=\int_{t}^{t+1}(-\dot{u}(\tau)(t+1-\tau)^{n+1}+u(\tau)(n+1)(t+1-\tau)^{n})d\tau.
\end{eqnarray*}
Then one has
\begin{eqnarray}
|u(t)|\leq\frac{1}{\sqrt{2n+3}}\left(\int_{t}^{t+1}|\dot{u}(\tau)|^{2}d\tau\right)^{\frac{1}{2}}+(n+1)\left(\int_{t}^{t+1}|u(\tau)|d\tau\right).\label{20}
\end{eqnarray}
For any $u$, $v\in K$ with $|t|\geq R+1$, where $R\geq r_{0}$, it
follows from $(L1)$, $(L2)$ and (\ref{20}) that
\begin{eqnarray*}
&&|u(t)-v(t)|\nonumber\\
&\leq&\frac{1}{\sqrt{2n+3}}\left(\int_{\Gamma_{R}^{c}}|\dot{u}(\tau)-\dot{v}(\tau)|^{2}d\tau\right)^{\frac{1}{2}}+(n+1)\int_{[t,t+1]\bigcap D_{M}(R)}|u(\tau)-v(\tau)|d\tau\nonumber\\
&&+(n+1)\int_{[t,t+1]\setminus D_{M}(R)}|u(\tau)-v(\tau)|d\tau\\
&\leq&\frac{\sqrt{2}}{\sqrt{2n+3}}\left(\|u\|_{o}+\|v\|_{o}\right)+(n+1)\int_{D_{M}(R)}|u(\tau)-v(\tau)|d\tau\nonumber\\
&&+(n+1)\sqrt{\mbox{meas}([t,t+1]\setminus D_{M}(R))}\left(\int_{[t,t+1]\setminus D_{M}(R)}|u(\tau)-v(\tau)|^{2}d\tau\right)^{\frac{1}{2}}\nonumber\\
&\leq&\frac{2\sqrt{2}M_{1}}{\sqrt{2n+3}}+(n+1)\sqrt{\mbox{meas}(D_{M}(R))}\left(\int_{D_{M}(R)}|u(\tau)-v(\tau)|^{2}d\tau\right)^{\frac{1}{2}}\nonumber\\
&&+(n+1)\left(\int_{\Gamma_{R}^{c}\setminus D_{M}(R)}|u(\tau)-v(\tau)|^{2}d\tau\right)^{\frac{1}{2}}\nonumber\\
&\leq&\frac{2\sqrt{2}M_{1}}{\sqrt{2n+3}}+(n+1)\sqrt{\mbox{meas}(D_{M}(R))}\left(\int_{\mathbb{R}}|u(\tau)-v(\tau)|^{2}d\tau\right)^{\frac{1}{2}}\nonumber\\
&&+\frac{n+1}{\sqrt{M}}\left(\int_{\Gamma_{R}^{c}\setminus D_{M}(R)}|\tau|^{-1}\ln^{-\alpha}|\tau|(L(\tau)(u(\tau)-v(\tau)),u(\tau)-v(\tau))d\tau\right)^{\frac{1}{2}}\nonumber\\
&\leq&\frac{2\sqrt{2}M_{1}}{\sqrt{2n+3}}+2(n+1)\sqrt{\mbox{meas}(D_{M}(R))}\left(\|u\|_{o}+\|v\|_{o}\right)\nonumber\\
&&+\frac{n+1}{\sqrt{M}R^{\frac{1}{2}}\ln^{\frac{\alpha}{2}}R}\left(\int_{\Gamma_{R}^{c}\setminus D_{M}(R)}(L(\tau)(u(\tau)-v(\tau)),u(\tau)-v(\tau))d\tau\right)^{\frac{1}{2}}\nonumber\\
&\leq&\frac{2\sqrt{2}M_{1}}{\sqrt{2n+3}}+4M_{1}(n+1)\sqrt{\mbox{meas}(D_{M}(R))}+\frac{(n+1)}{\sqrt{M}R^{\frac{1}{2}}\ln^{\frac{\alpha}{2}}R}\left(\|u\|_{o}+\|v\|_{o}\right)\nonumber\\
&\leq&\frac{2\sqrt{2}M_{1}}{\sqrt{2n+3}}+4M_{1}(n+1)\sqrt{\mbox{meas}(D_{M}(R))}+\frac{2M_{1}(n+1)}{\sqrt{M}R^{\frac{1}{2}}\ln^{\frac{\alpha}{2}}R}.
\end{eqnarray*}
Then there exists $n_{0}$ such that
$\frac{2\sqrt{2}M_{1}}{\sqrt{2n+3}}<\frac{1}{3}\varepsilon$ for any
$n\geq n_{0}$. Then we can choose $R$ large enough such that the last
two terms in above inequality are both smaller than
$\frac{1}{3}\varepsilon$. Then one can deduce that
\begin{eqnarray*}
\max_{t\in \Gamma_{R+1}^{c}}|u(t)-v(t)|\leq\varepsilon.
\end{eqnarray*}
Similar to the case $p=1$,  $K$ is precompact in
$L^{\infty}(\mathbb{R},\mathbb{R}^{N})$.

Finally, we will show our conclusion holds without $l(t)\geq1$. From
$(L1)$, we can see $L(t)$ is bounded from below. Hence there exists a
constant $a>0$ such that $l(t)+a\geq1$ for all $t\in \mathbb{R}$. On
$\mathfrak{D}((\mathcal{A}+a)^{\frac{1}{2}})$, we introduce a norm
\begin{eqnarray*}
\|u\|_{\sim}^{2}=\left\|(\mathcal{A}+a)^{\frac{1}{2}}u\right\|_{2}^{2}+\|u\|_{2}^{2}.
\end{eqnarray*}
Then it follows from the proof of Lemma 2.1 in \cite{2} that
$\|\cdot\|_{\sim}$ is equivalent to $\|\cdot\|_{o}$ and $\|\cdot\|$.
 We obtain our conclusion.\hfill$\Box$

\begin{lemma.}\label{le2.0}
 Under the conditions of Theorem \ref{th1.1}, we see $I$ is of $C^{1}$ class and
\begin{eqnarray*}
\langle I'(u),u \rangle=\|u^{+}\|^{2}-\|u^{-}\|^{2}-\int_{\mathbb{R}}(\nabla
W(t,u(t)),u(t)))dt,\ \ \mbox{ $\forall u$ $\in E$}.
\end{eqnarray*}
\end{lemma.}

{\bf Proof.} The proof is similar to Lemma 2.2 in
\cite{16}.\hfill$\Box$

Next, we show problem (\ref{1}) has at least a nontrivial homoclinic
solution by the Generalized Mountain Pass Theorem under $(PS)$
condition (see Theorem 5.29 in \cite{6}). First, we show $I$ satisfies
the $(PS)$ condition.

\begin{lemma.}\label{le2.2} Suppose $(L1)$-$(L2)$ and $(W1)$-$(W6)$ hold, then any $(PS)$ sequence of $I$ is bounded.
\end{lemma.}

{\bf Proof.} Assume that $\{u_{n}\} \subset E$ is a $(PS)$ sequence,
then there is a $M_{3} > 0$ such that
\begin{eqnarray}
|I(u_{n})| \leq M_{3},\ \ \ \ \ \|I'(u_{n})\| \leq
M_{3}.\label{9}
\end{eqnarray}

Now we show the boundedness of $\{u_{n}\}$. We adopt an indirect
argument. Assuming that $\|u_{n}\|\rightarrow+\infty$ as $n\rightarrow
\infty$, then one can deduce from (\ref{9}), $(W2)$ and $(W3)$ that
\begin{eqnarray}
2M_{3}+M_{3}\|u_{n}\|&\geq&\langle I'(u_{n}),u_{n}\rangle-2I(u_{n})\nonumber\\
&=&\int_{\mathbb{R}}\widetilde{W}(t,u_{n}(t))dt
\nonumber\\
&=&\int_{\{t\in \mathbb{R}: |u_{n}|\leq r_{1}\}}\widetilde{W}(t,u_{n}(t))dt+\int_{\{t\in \mathbb{R}: |u_{n}|\geq r_{1}\}}\widetilde{W}(t,u_{n}(t))dt\nonumber\\
&\geq&-b_{3}\int_{\{t\in \mathbb{R}: |u_{n}|\leq r_{1}\}}|u_{n}(t)|dt+b_{2}\int_{\{t\in \mathbb{R}: |u_{n}|\geq r_{1}\}}|u_{n}(t)|^{\mu}dt.\label{8}
\end{eqnarray}
Let
\begin{eqnarray*}
v_{n}(t)=\left\{
\begin{array}{ll}
u_{n}(t)&\mbox{for $|u_{n}(t)|\leq r_{1}$}\\
0&\mbox{for $|u_{n}(t)|\geq r_{1}$}
\end{array}
\right.
\end{eqnarray*}
and
\begin{eqnarray*}
w_{n}(t)=u_{n}(t)-v_{n}(t).
\end{eqnarray*}
Then we have
\begin{eqnarray*}
b_{2}\|w_{n}\|_{\mu}^{\mu}\leq2M_{3}+M_{3}\|u_{n}\|+b_{3}\|u_{n}\|_{1}.
\end{eqnarray*}
Since $E^{-}\oplus E^{0}$ is of finite dimension, there is $M_{4}>0$
such that
\begin{eqnarray*}
\|u_{n}^{-}+u_{n}^{0}\|_{2}^{2}&=&(u_{n}^{-}+u_{n}^{0},u_{n})_{2}\nonumber\\
&=&(u_{n}^{-}+u_{n}^{0},v_{n})_{2}+(u_{n}^{-}+u_{n}^{0},w_{n})_{2}\nonumber\\
&\leq&\|u_{n}^{-}+u_{n}^{0}\|_{1}\|v_{n}\|_{\infty}+\|u_{n}^{-}+u_{n}^{0}\|_{L^{\frac{\mu}{\mu-1}}}\|w_{n}\|_{\mu}\nonumber\\
&\leq&M_{4}\|u_{n}^{-}+u_{n}^{0}\|_{2}(1+\|w_{n}\|_{\mu}),
\end{eqnarray*}
Moreover, there exist $M_{5}$, $M_{6}>0$ such that
\begin{eqnarray*}
\|u_{n}^{-}+u_{n}^{0}\|\leq M_{5}\|u_{n}^{-}+u_{n}^{0}\|_{2}\leq M_{4}M_{5}(1+\|w_{n}\|_{\mu})\leq M_{6}(1+\|u_{n}\|^{\frac{1}{\mu}}).
\end{eqnarray*}
Then one has
\begin{eqnarray}
\frac{\|u_{n}^{-}+u_{n}^{0}\|}{\|u_{n}\|}\rightarrow0\ \ \ \mbox{as}\ \ \ n\rightarrow\infty.\label{2}
\end{eqnarray}
Choosing
$0<s<\min\left\{1,\frac{\mu}{2},\frac{\mu}{2(\mu-1)}\right\}$, we can
deduce from (\ref{4}) and (\ref{8}) that
\begin{eqnarray}
\int_{\{t\in \mathbb{R}: |u_{n}|\geq r_{1}\}}|u_{n}(t)|^{2}dt&=&\int_{\{t\in \mathbb{R}: |u_{n}|\geq r_{1}\}}|u_{n}(t)|^{2s}|u_{n}(t)|^{2(1-s)}dt\nonumber\\
&\leq&\left(\int_{\{t\in \mathbb{R}: |u_{n}|\geq r_{1}\}}|u_{n}(t)|^{\mu}dt\right)^{\frac{2s}{\mu}}\left(\int_{\{t\in \mathbb{R}: |u_{n}|\geq r_{1}\}}|u_{n}(t)|^{\frac{2\mu(1-s)}{\mu-2s}}dt\right)^{\frac{\mu-2s}{\mu}}\nonumber\\
&\leq&M_{7}\left(1+\int_{\{t\in \mathbb{R}: |u_{n}|\leq r_{1}\}}|u_{n}(t)|dt+\|u_{n}\|\right)^{\frac{2s}{\mu}}\|u_{n}\|_{\frac{2\mu(1-s)}{\mu-2s}}^{2(1-s)}\nonumber\\
&\leq&M_{8}\left(1+\|u_{n}\|\right)^{\frac{2s}{\mu}}\|u_{n}\|^{2(1-s)}\nonumber\\
&\leq&M_{9}\left(\|u_{n}\|^{2(1-s)}+\|u_{n}\|^{\frac{2(s+\mu-s\mu)}{\mu}}\right)\label{15}
\end{eqnarray}
for some $M_{7}$, $M_{8}$, $M_{9}>0$. By (\ref{9}), we obtain
\begin{eqnarray}
M_{3}(1+\|u_{n}\|)&\geq&\|I'(u_{n})\|(1+\|u_{n}\|)\nonumber\\
&\geq&\|I'(u_{n})\|(1+\|u_{n}^{+}\|)\nonumber\\
&\geq&|\langle I'(u_{n}),u_{n}^{+}\rangle|\nonumber\\
&=&\left|\|u_{n}^{+}\|^{2}-\int_{\mathbb{R}}(\nabla W(t,u_{n}(t)),u_{n}^{+}(t))dt\right|\nonumber\\
&\geq&\|u_{n}^{+}\|^{2}-\int_{\mathbb{R}}|\nabla W(t,u_{n}(t)||u_{n}^{+}(t)|dt.\label{3}
\end{eqnarray}
By $(W1)$, (\ref{15}) and (\ref{3}), we deduce
\begin{eqnarray*}
&&\|u_{n}^{+}\|^{2}\nonumber\\
&\leq&\int_{\mathbb{R}}|\nabla W(t,u_{n}(t)||u_{n}^{+}(t)|dt+M_{3}(1+\|u_{n}\|)\nonumber\\
&\leq&b_{1}\int_{\mathbb{R}}|u_{n}(t)||u_{n}^{+}(t)|dt+M_{3}(1+\|u_{n}\|)\nonumber\\
&=&b_{1}\left(\int_{\{t\in \mathbb{R}: |u_{n}|\leq r_{1}\}}|u_{n}(t)||u_{n}^{+}(t)|dt+\int_{\{t\in \mathbb{R}: |u_{n}|\geq r_{1}\}}|u_{n}(t)||u_{n}^{+}(t)|dt\right)+M_{3}(1+\|u_{n}\|)\nonumber\\
&\leq&b_{1}\left(r_{1}\int_{\{t\in \mathbb{R}: |u_{n}|\leq r_{1}\}}|u_{n}^{+}(t)|dt+\left(\int_{\{t\in \mathbb{R}: |u_{n}|\geq r_{1}\}}|u_{n}(t)|^{2}dt\right)^{\frac{1}{2}}\left(\int_{\{t\in \mathbb{R}: |u_{n}|\geq r_{1}\}}|u_{n}^{+}(t)|^{2}dt\right)^{\frac{1}{2}}\right)\nonumber\\
&&+M_{3}(1+\|u_{n}\|)\nonumber\\
&\leq&M_{10}\left(\|u_{n}^{+}\|+\|u_{n}^{+}\|\left(\|u_{n}\|^{1-s}+\|u_{n}\|^{\frac{s+\mu-s\mu}{\mu}}\right)\right)+M_{3}(1+\|u_{n}\|)\nonumber\\
&\leq&M_{11}\left(\|u_{n}\|+\|u_{n}\|^{2-s}+\|u_{n}\|^{\frac{s+2\mu-s\mu}{\mu}}\right)+M_{3}(1+\|u_{n}\|)
\end{eqnarray*}
for some $M_{10}$, $M_{11}>0$. Hence one obtains
\begin{eqnarray}
\frac{\|u_{n}^{+}\|^{2}}{\|u_{n}\|^{2}}\rightarrow0\ \ \ \mbox{as}\ \ \ n\rightarrow\infty.\label{5}
\end{eqnarray}
Formula (\ref{2}) and (\ref{5}) show a contradiction. Therefore
$\{u_{n}\}$ is bounded.\hfill$\Box$

\begin{lemma.}\label{le2.3} Suppose $(L1)$-$(L2)$ and $(W1)$-$(W6)$ hold, then we can find $\alpha$, $\varrho>0$
such that $I\mid_{S} \geq \alpha$, where $S$ is a sphere in $E^{+}$ of
radius $\varrho$, i.e. $S=\{u\in E^{+}|\ \|u\|=\varrho\}$.
\end{lemma.}

{\bf Proof.} Set
\begin{eqnarray*}
\ \varrho =\frac{r_{2}}{C_{\infty}},\ \ \ \ \ \ \ \
\alpha=\frac{\sigma_{0}}{2\lambda_{\bar{n}+1}}\varrho^{2}>0.\label{14}
\end{eqnarray*}
One sees that $0 < \|u\|_{\infty} \leq r_{2}$, $\forall u\in S$. By
(\ref{7}) and $(W4)$, one can see that
\begin{eqnarray}
\ I(u)&=&
\frac{1}{2}\|u\|^{2}-\int_{\mathbb{R}}W(t,u(t))dt{\nonumber}\\
&\geq& \mbox{}
\frac{1}{2}\|u\|^{2}-\frac{\lambda_{\bar{n}+1}-\sigma_{0}}{2}\int_{\mathbb{R}}|u(t)|^{2}dt{\nonumber}\\
&\geq& \mbox{}
\frac{1}{2}\left(1-\frac{\lambda_{\bar{n}+1}-\sigma_{0}}{\lambda_{\bar{n}+1}}\right)\|u\|^{2}{\nonumber}\\
&=& \mbox{} \frac{\sigma_{0}}{2\lambda_{\bar{n}+1}}\|u\|^{2}.\label{24}
\end{eqnarray}
By the definitions of $\varrho$ and $\alpha$, ~(\ref{24}) implies
$I\mid_{S} \geq \alpha$. We obtain our
conclusion.~~~~~~~~~~~~~~~~~~~~~~~~~~~~~~~~~~~$\Box$

\begin{lemma.}\label{le2.4} Suppose $(L1)$-$(L2)$ and $(W1)$-$(W6)$ hold, then there exists $\rho>\varrho>0$ such that
 $\sup_{\partial Q} I(u) \leq 0$, where
\begin{eqnarray*}
Q=\{\zeta e_{\bar{n}+1}|\
0\leq\zeta\leq\rho\}\oplus\{u\in E^{-}\oplus E^{0}|\ \|u\|\leq
\rho\}.
\end{eqnarray*}
\end{lemma.}

{\bf Proof.} It follows from $(W5)$ and $(W6)$ that
\begin{eqnarray}
W(t,x)\geq \left(\frac{\lambda_{\bar{n}+1}}{2}+\varepsilon_{0}\right)(|x|^{2}-r_{\infty}|x|)\label{11}
\end{eqnarray}
for $(t,x)\in\mathbb{R}\times\mathbb{R}^{N}$. By the definition of
$Q$, for any $u\in Q$, we have
\begin{eqnarray}
u=\zeta e_{\bar{n}+1}+v\label{6}
\end{eqnarray}
for some $0\leq\zeta\leq\rho$ and $v\in E^{-}\oplus E^{0}$ with
$\|v\|\leq \rho$. It can be deduced from (\ref{7}), (\ref{11})
and~(\ref{6}) that
\begin{eqnarray}
I(u)
&\leq& \frac{1}{2}\|u^{+}\|^{2}-\left(\frac{\lambda_{\bar{n}+1}}{2}+\varepsilon_{0}\right)\int_{\mathbb{R}}|u(t)|^{2}dt+M_{12}\int_{\mathbb{R}}|u(t)|dt{\nonumber}\\
&\leq& \frac{\zeta^{2}}{2}\|e_{\bar{n}+1}\|^{2}-\left(\frac{\lambda_{\bar{n}+1}}{2}+\varepsilon_{0}\right)\int_{\mathbb{R}}(|\zeta e_{\bar{n}+1}(t)|^{2}+|v(t)|^{2})dt+M_{13}\|u\|{\nonumber}\\
&\leq& \frac{\zeta^{2}}{2}\|e_{\bar{n}+1}\|^{2}-\zeta^{2}\left(\frac{\lambda_{\bar{n}+1}}{2}+\varepsilon_{0}\right)\|e_{\bar{n}+1}\|_{2}^{2}-\left(\frac{\lambda_{\bar{n}+1}}{2}+\varepsilon_{0}\right)\|v\|_{2}^{2}+M_{13}\|u\|{\nonumber}\\
&\leq& -\frac{\zeta^{2}\varepsilon_{0}}{\lambda_{\bar{n}+1}}\|e_{\bar{n}+1}\|^{2}-\left(\frac{\lambda_{\bar{n}+1}}{2}+\varepsilon_{0}\right)\|v\|_{2}^{2}+M_{13}\|u\|.\label{10}
\end{eqnarray}
for some $M_{12}$, $M_{13}>0$ and any $u\in Q$. By the definition of
$Q$, we shows
\begin{eqnarray*}
\partial Q=Q_{1}\cup Q_{2}\cup Q_{3},
\end{eqnarray*}
where
\begin{eqnarray*}
\begin{array}{ll}
Q_{1}=\{u\in E^{-}\oplus E^{0}|\ \|u\|\leq
\rho\}.\\
Q_{2}=\rho e_{\bar{n}+1}\oplus\{u\in E^{-}\oplus E^{0}|\ \|u\|\leq
\rho\}.\\
Q_{3}=\{\zeta e_{\bar{n}+1}|\
0\leq\zeta\leq\rho\}\oplus\{u\in E^{-}\oplus E^{0}|\ \|u\|=
\rho\}.
\end{array}
\end{eqnarray*}
Therefore, we deduce from $(W6)$ that $\sup_{u\in Q_{1}} I(u) \leq 0$.
By (\ref{10}), one has
\begin{eqnarray*}
\sup_{u\in Q_{2}} I(u)
\leq -\frac{\rho^{2}\varepsilon_{0}}{\lambda_{\bar{n}+1}}\|e_{\bar{n}+1}\|^{2}+M_{13}\rho,
\end{eqnarray*}
which implies that $\sup_{u\in Q_{2}} I(u)\leq0$ for $\rho$ large
enough. Finally, one can deduce from (\ref{10}) that
\begin{eqnarray*}
\sup_{u\in Q_{3}} I(u)
\leq -M_{14}\rho^{2}+M_{13}\rho
\end{eqnarray*}
for some $M_{14}>0$ since $\dim(E^{-}\oplus E^{0})<\infty$, which
implies that $\sup_{u\in Q_{3}} I(u)\leq0$ for $\rho$ large enough.
Then we obtain our
conclusion.~~~~~~~~~~~~~~~~~~~~~~~~~~~~~~~~~~~~~~~~~~~~~~~~~~~~~~~~~~~~~~~~~~~~~~~~~~~~~~~~~~~~~~~~~~~~~~$\Box$

\vspace{0.3cm}{\bf Proof of Theorem \ref{th1.1}.} From Lemmas
\ref{le2.3}, \ref{le2.4} and the Generalized Mountain Pass theorem,
there exists a sequence $\{u_{k}\}$ such that that $\{I(u_{k})\}$ is
bounded, $I'(u_{k})\rightarrow 0$ as $k \rightarrow \infty$. It can be
obtained from Lemma \ref{le2.2} that $\{u_{k}\}$ is a bounded sequence
in $E$. Then we can find a subsequence, which is still denoted by
$\{u_{k}\}$, such that $u_{k}\rightharpoonup u_{0}$ in $E$, which
implies that $u_{k}\rightarrow u_{0}$ in
$L^{2}(\mathbb{R},\mathbb{R}^{N})$. One can imply that
\begin{eqnarray*}
&&\|u_{k}^{+}-u_{0}^{+}\|^{2}\nonumber\\&=&\langle I'(u_{k})-I'(u_{0}),u_{k}^{+}-u_{0}^{+} \rangle+\int_{\mathbb{R}}(\nabla
W(t,u_{k}(t))-\nabla
W(t,u_{0}(t)),u_{k}^{+}(t)-u_{0}^{+}(t))dt\nonumber\\
&\leq&\|I'(u_{k})\|(\|u_{k}^{+}\|+\|u_{0}^{+}\|)-\langle I'(u_{0}),u_{k}^{+}-u_{0}^{+} \rangle\nonumber\\
&&+\int_{\mathbb{R}}|\nabla
W(t,u_{k}(t))-\nabla
W(t,u_{0}(t))||u_{k}^{+}(t)-u_{0}^{+}(t)|dt\nonumber\\
&\leq&o(1)+b_{1}\int_{\mathbb{R}}(|u_{k}(t)|+|u_{0}(t)|)|u_{k}^{+}(t)-u_{0}^{+}(t)|dt\nonumber\\
&\leq&o(1)+2b_{1}(\|u_{k}\|_{2}+\|u_{0}\|_{2})\left(\int_{\mathbb{R}}|u_{k}^{+}(t)-u_{0}^{+}(t)|^{2}dt\right)^{\frac{1}{2}}\nonumber\\
&\leq&o(1)+2b_{1}C_{2}(\|u_{k}\|+\|u_{0}\|)\left(\int_{\mathbb{R}}|u_{k}^{+}(t)-u_{0}^{+}(t)|^{2}dt\right)^{\frac{1}{2}}\nonumber\\
&\rightarrow&0\ \ \ \mbox{as}\ \ \ k\rightarrow\infty.
\end{eqnarray*}
Then $u_{k}\rightarrow u_{0}$ as $n\rightarrow\infty$ since
$E^{-}\oplus E^{0}$ is of finite dimensions. Therefore, we have
$I'(u_{0})=0$. Finally, we need to prove that $u_{0}\in
\mathfrak{D}(\mathcal{H})$. It is easy to see that $u_{0}\in C^{2}$
satisfies (\ref{1}). Then by (\ref{1}), (\ref{4}) and $(W_{1})$, we
have
\begin{eqnarray*}
\|\mathcal{H} u_{0}\|_{2}^{2}=\int_{\mathbb{R}}|\nabla
W(t,u_{0}(t))|^{2}dt\leq b_{1}^{2}\int_{\mathbb{R}}|u_{0}(t)|^{2}dt\leq b_{1}^{2}C_{2}^{2}\|u_{0}\|^{2}<+\infty.
\end{eqnarray*}
Thus $u_{0}\in \mathfrak{D}(\mathcal{H})$. We obtain our conclusion.
~~~~~~~~~~~~~~~~~~~~~~~~~~~~~~~~~~~~~~~~~~~~~~~~~~~~~~~~~~~~~~~~~~~~~~~~~~~~~~$\Box$

\section{Proof of Theorem \ref{th1.2}}

Subsequently, it will be shown that $I$ possesses infinitely many
critical points if $W(t,x)$ is even in $x$, which are gotten by the
Variant Fountain Theorem by Zou in \cite{30}. For any $j\in
\mathbb{N}$, define $X_{j}=span\{e_{j}\}$. Set
$Y_{k}=\bigoplus_{j=1}^{k}X_{j}$, $Z_{k}=\overline{\bigoplus_{j=
k}^{\infty}X_{j}}$, $B_{k}$ be a ball in $Y_{k}$ with radius
$\rho_{k}$ and centered at origin. Let $N_{k}$ be the boundary of a
ball in $Z_{k}$ with radius $r_{k}$ and centered at origin for
$\rho_{k}>r_{k}>0$. Obviously, we have $E=\overline{\bigoplus_{j\in
\mathbb{N}}X_{j}}$. Next, we show that $I$ satisfies the geometrical
structure of the Variant Fountain Theorem. Set
\begin{eqnarray*}
T(u)=\frac{1}{2}\|u^{+}\|^{2}, \ \
\ P(u)=\frac{1}{2}\|u^{-}\|^{2}+\int_{\mathbb{R}}W(t,u(t))dt
\end{eqnarray*}
and
\begin{eqnarray*}
I_{\lambda}(u)=T(u)-\lambda P(u).
\end{eqnarray*}

Suppose that $I_{\lambda}(u)$ satisfies

$(A_{1})$ $I_{\lambda}$ is a bounded map uniformly in
$\lambda\in[1,2]$. And $I_{\lambda}(-u)=I_{\lambda}(u)$ for
$\forall(\lambda,u)\in[1,2]\times E$;

$(A_{2})$ $P(u)$ is nonnegative on $E$; $T(u)$ or $P(u)$ is coercive;
or

$(A_{2}')$ $P(u)$ is nonpositive on $E$; $-P(u)$ is coercive.

For $k>2$, let
$$\Sigma_{k}:=\{\pi\in\mathcal{C}(B_{k},E): \pi\ \mbox{is an odd map},\ \pi|_{\partial B_{k}}=id\},$$
$$c_{k}(\lambda):=\inf_{\pi\in \Sigma_{k}}\max_{u\in
B_{k}}I_{\lambda}(\pi(u)),$$
$$b_{k}(\lambda):=\inf_{u\in Z_{k}, \|u\|=r_{k}}I_{\lambda}(u),$$
$$a_{k}(\lambda):=\max_{u\in Y_{k}, \|u\|=\rho_{k}}I_{\lambda}(u),$$

\begin{lemma.}(Zou\cite{30})\label{le4.1} Assume $(A_{1})$ and $(A_{2})$(or
$(A_{2}')$). If $b_{k}(\lambda)>a_{k}(\lambda)$ for $\forall\lambda\in
[1,2]$, then $c_{k}(\lambda)\geq b_{k}(\lambda)$ for
$\forall\lambda\in [1,2]$. Furthermore, for a.e. $\lambda\in [1,2]$,
there is $\{u_{n}^{k}(\lambda)\}_{n=1}^{\infty}$ such that
$$\sup_{n}\|u_{n}^{k}(\lambda)\|<\infty,\ I_{\lambda}'(u_{n}^{k}(\lambda))\rightarrow0\ \ \mbox{and}\ \ I_{\lambda}(u_{n}^{k}(\lambda))\rightarrow c_{k}(\lambda)\ \ \mbox{as}\ \ n\rightarrow\infty.$$
\end{lemma.}

\begin{lemma.}\label{le4.3} Suppose $(L1)$-$(L2)$, $(W1)$-$(W6)$ hold,
then we can find $r_{k}>0$, $\overline{b}_{k}\rightarrow+\infty$ as
$k\rightarrow\infty$ such that $b_{k}(\lambda)\geq\overline{b}_{k}$
for $\forall\lambda\in[1,2]$.
\end{lemma.}

{\bf Proof.} For $k\in \mathbb{N}$ and $p\in[1,+\infty]$, we set
\begin{eqnarray*}
\beta_{k}(p)=\sup_{u\in Z_{k},\ \|u\|=1}\|u\|_{L^{p}}.
\end{eqnarray*}
Obviously, $\beta_{k}(p)\rightarrow0^{+}$ as $k\rightarrow\infty$.
Then, one has $Z_{k}\subset E^{+}$ and
$\beta_{k}(2)\leq\frac{1}{4(\lambda_{\bar{n}+1}-\sigma_{0})}$ for $k$
large enough. By $(W1)$ and $(W4)$, we conclude that
\begin{eqnarray}
|W(t,x)|\leq \frac{\lambda_{\bar{n}+1}-\sigma_{0}}{2}|x|^{2}+M_{15}|x|^{q}\ \ \ \mbox{for}\ \ \forall (t,x)\in \mathbb{R}\times \mathbb{R}^{N}\label{22}
\end{eqnarray}
for some $M_{15}>0$ and $q>2$. For any $u\in Z_{k}$, we can deduce
that
\begin{eqnarray*}
I_{\lambda}(u)&=&
\frac{1}{2}\|u\|^{2}-\lambda\int_{\mathbb{R}^{N}}W(t,u(t))dt{\nonumber}\\
&\geq& \mbox{}
\frac{1}{2}\|u\|^{2}-(\lambda_{\bar{n}+1}-\sigma_{0})\int_{\mathbb{R}^{N}}|u(t)|^{2}dt-2M_{15}\int_{\mathbb{R}^{N}}|u(t)|^{q}dt{\nonumber}\\
&\geq& \mbox{}
\frac{1}{2}\left(1-2(\lambda_{\bar{n}+1}-\sigma_{0})\beta_{k}(2)\right)\|u\|^{2}-2M_{15}\beta_{k}^{q}(q)\|u\|^{q}{\nonumber}\\
&=& \mbox{}
\frac{1}{4}\|u\|^{2}-2M_{15}\beta_{k}^{q}(q)\|u\|^{q}.\label{96}
\end{eqnarray*}
Let
\begin{eqnarray*}
r_{k}=\left(\frac{1}{16M_{15}\beta_{k}^{q}(q)}\right)^{\frac{1}{q-2}}.
\end{eqnarray*}
Obviously, we have $r_{k}\rightarrow+\infty$ as $k\rightarrow\infty$.
When $k$ large enough, we can conclude
\begin{eqnarray*}
b_{k}(\lambda)&=&\inf_{u\in
Z_{k},\|u\|=r_{k}}I_{\lambda}(u)\nonumber\\
&\geq& \mbox{}\frac{1}{4}r_{k}^{2}-2M_{15}\beta_{k}^{q}(q)r_{k}^{q}\nonumber\\
&\geq& \mbox{}r_{k}^{2}\left(\frac{1}{4}-2M_{15}\beta_{k}^{q}(q)r_{k}^{q-2}\right)\nonumber\\
&=& \mbox{}\frac{1}{8}r_{k}^{2}\nonumber\\
&\doteq&\overline{b}_{k}.
\end{eqnarray*}
Then $\overline{b}_{k}\rightarrow+\infty$ as
$k\rightarrow\infty$.~~~~~~~~~~~~~~~~~~~~~~~~~~~~~~~~~~~~~~~~~~~~~~~~~~~~~~~~~~~~~~~~~~~~~~~~~~~~~~~~~~~~~~~~~~~~~~~~~~~$\Box$

\begin{lemma.}\label{le4.4} We can find $\rho_{k}$ big enough such
that $a_{k}(\lambda)\leq0$ for $\forall\lambda\in[1,2]$.
\end{lemma.}

{\bf Proof.} For $\forall u\in Y_{k}\setminus\{0\}$, $\xi\in(0,1)$ and
$\delta>0$, set
\begin{eqnarray*}
\Gamma_{\delta}(u)=\{t\in \mathbb{R}:\ |u|^{\xi}\geq\delta \|u\|\}.
\end{eqnarray*}
Subsequently, we prove that there exists $\varepsilon_{1}>0$ such that
\begin{eqnarray}
\mbox{meas}(\Gamma_{\varepsilon_{1}}(u))\geq\varepsilon_{1}\ \ \ \mbox{for all}\ \ u \in Y_{k}\setminus\{0\}.\label{60}
\end{eqnarray}
Otherwise, there exists $\{u_{n}\}_{n\in \mathbb{N}}$ such that
\begin{eqnarray}
\mbox{meas}(\Gamma_{\frac{1}{n}}(u_{n}))\leq \frac{1}{n}.\label{58}
\end{eqnarray}
Set $v_{n}=\frac{|u_{n}|}{\|u_{n}\|}$, then $\|v_{n}\|=1$ for $\forall
n\in \mathbb{N}$, which implies that there exists $v_{0}\in Y_{k}$
such that $\|v_{0}\|=1$ and $v_{n}\rightarrow v_{0}$ in
$L^{2}(\mathbb{R},\mathbb{R}^{N})$. Then one obtains
\begin{eqnarray}
\int_{\mathbb{R}}|v_{n}(t)-v_{0}(t)|^{2}dt\rightarrow0\ \ \ \mbox{as}\ \
n\rightarrow\infty.\label{59}
\end{eqnarray}
Furthermore, we can obtain that there exist $\tau_{1}$, $\tau_{2}>0$
such that
\begin{eqnarray}
\mbox{meas}(\Gamma_{\tau_{1}}(v_{0}))\geq\tau_{2}.\label{57}
\end{eqnarray}
If not, we have $\mbox{meas}(\Gamma_{\frac{1}{n}}(v_{0}))=0$ for
$\forall n\in \mathbb{N}$. Then one gets
\begin{eqnarray*}
0\leq\int_{\mathbb{R}}|v_{0}(t)|^{1+\xi}dt
=\int_{\mathbb{R}}|v_{0}(t)|^{\xi}|v_{0}(t)|dt
\leq\frac{1}{n}\int_{\mathbb{R}}|v_{0}(t)|dt
\leq\frac{C_{1}}{n}\rightarrow0\ \ \ \mbox{as}\ \ n\rightarrow\infty,
\end{eqnarray*}
which contradicts $\|v_{0}\|=1$. Then (\ref{57}) holds. For $n$ large
enough, we can deduce that
\begin{eqnarray}
|v_{n}(t)-v_{0}(t)|^{2}\geq||v_{n}(t)|-|v_{0}(t)||^{2}\geq\left(\tau_{1}^{\frac{1}{\xi}}-\left(\frac{1}{n}\right)^{\frac{1}{\xi}}\right)^{2}\geq\frac{1}{4}\tau_{1}^{\frac{2}{\xi}}
\end{eqnarray}
for all
$t\in\Gamma_{\frac{1}{n}}^{c}(v_{n})\bigcap\Gamma_{\tau_{1}}(v_{0})$.
By (\ref{58}) and (\ref{57}), we obtain that
\begin{eqnarray}
m\left(\Gamma_{\frac{1}{n}}^{c}(v_{n})\bigcap\Gamma_{\tau_{1}}(v_{0})\right)&=&m\left(\Gamma_{\tau_{1}}(v_{0})\setminus\left(\Gamma_{\frac{1}{n}}(v_{n})\bigcap\Gamma_{\tau_{1}}(v_{0})\right)\right){\nonumber}\\
&\geq&\mbox{}m\left(\Gamma_{\tau_{1}}(v_{0})\right)-m\left(\Gamma_{\frac{1}{n}}(v_{n})\bigcap\Gamma_{\tau_{1}}(v_{0})\right){\nonumber}\\
&\geq&\mbox{}\tau_{2}-\frac{1}{n}.
\end{eqnarray}
Consequently, we have
\begin{eqnarray}
\int_{\mathbb{R}}|v_{n}(t)-v_{0}(t)|^{2}dt&\geq&\int_{\Gamma_{\frac{1}{n}}^{c}(v_{n})\bigcap\Gamma_{\tau_{1}}(v_{0})}|v_{n}(t)-v_{0}(t)|^{2}dt{\nonumber}\\
&\geq&\mbox{}\frac{1}{4}\tau_{1}^{\frac{2}{\xi}}\left(\tau_{2}-\frac{1}{n}\right){\nonumber}\\
&\geq&\mbox{}\frac{1}{8}\tau_{1}^{\frac{2}{\xi}}\tau_{2}{\nonumber}\\
&>&\mbox{}0
\end{eqnarray}
for $n$ big enough, which contradicts (\ref{59}). Therefore (\ref{60})
holds.

For any $u\in Y_{k}$ and $t\in \Gamma_{\varepsilon_{1}}(u)$ with
$\|u\|\geq \frac{r_{\infty}^{\xi}}{\varepsilon_{1}}$, it follows from
$(W5)$ that
\begin{eqnarray}
W(t,u(t))\geq \left(\frac{\lambda_{\bar{n}+1}}{2}+\varepsilon_{0}\right)|u(t)|^{2}\geq
\left(\frac{\lambda_{\bar{n}+1}}{2}+\varepsilon_{0}\right)\varepsilon_{1}^{\frac{2}{\xi}}\|u\|^{\frac{2}{\xi}}.\label{12}
\end{eqnarray}
We can choose
$\rho_{k}>\max\left\{\frac{r_{\infty}^{\xi}}{\varepsilon_{1}},r_{k}\right\}$,
then for any $u\in Y_{k}$ with $\|u\|=\rho_{k}$, it follows from
(\ref{7}), $(W6)$, (\ref{60}) and (\ref{12}) that
\begin{eqnarray*}
I_{\lambda}(u)&=&
\frac{1}{2}\|u^{+}\|^{2}-\lambda\left(\frac{1}{2}\|u^{-}\|^{2}+\int_{\mathbb{R}}W(t,u(t))dt\right){\nonumber}\\
&\leq& \mbox{}\frac{1}{2}\|u\|^{2}-\int_{\Gamma_{\varepsilon_{1}}(u)}W(t,u(t))dt{\nonumber}\\
&\leq& \mbox{}\frac{1}{2}\|u\|^{2}-\left(\frac{\lambda_{\bar{n}+1}}{2}+\varepsilon_{0}\right)\int_{\Gamma_{\varepsilon_{1}}(u)}|u(t)|^{2}dx{\nonumber}\\
&\leq& \mbox{}\frac{1}{2}\|u\|^{2}-\left(\frac{\lambda_{\bar{n}+1}}{2}+\varepsilon_{0}\right)\varepsilon_{1}^{\frac{\xi+2}{\xi}}\|u\|^{\frac{2}{\xi}}{\nonumber}\\
&\leq& \mbox{}0,
\end{eqnarray*}
which means $a_{k}(\lambda)\leq0$ for $\rho_{k}$ large enough since
$\xi\in(0,1)$.
~~~~~~~~~~~~~~~~~~~~~~~~~~~~~~~~~~~~~~~~~~~~~~~~~~~~~~~~~~~~~~~~~~~~~~$\Box$

By Lemma 3.1, we can obtain

\begin{lemma.}\label{le4.5}
There exist $\lambda_{n}\rightarrow1$ as $n\rightarrow\infty$,
$\{u_{n}(k)\}_{n=1}^{\infty}\subset E$ sastisfying
$I_{\lambda_{n}}'(u_{n}(k))=0$ and
$I_{\lambda_{n}}(u_{n}(k))\in[\overline{b}_{k},\overline{c}_{k}]$,
 where $\overline{c}_{k}=\sup_{u\in B_{k}}I(u)$.
\end{lemma.}

{\bf Proof of Theorem \ref{th1.2}.} From Lemmas \ref{le2.2} and
\ref{le4.5}, $\{u_{n}(k)\}_{n=1}^{\infty}$ is a bounded sequence.
Similar to Theorem \ref{th1.1}, we can conclude that
$u_{n}(k)\rightarrow u(k)$ as $n\rightarrow\infty$, which implies that
$I_{1}$ possesses a critical point $u_{k}$ with
$I_{1}(u(k))\in[\overline{b}_{k},\overline{c}_{k}]$. Therefore,
 we obtain a sequence of critical points of $I_{1}$ since
 $\overline{b}_{k}\rightarrow+\infty$ as $k\rightarrow\infty$.\hfill$\Box$

\def\refname{References}

\end{document}